\documentclass{pspum-l}
\usepackage{amsmath,cite,epsf,epic}
\usepackage{eepic}

\copyrightinfo{2000}{American Mathematical Society}

\makeatletter
\renewcommand*\subjclass[2][2000]{%
  \def\@subjclass{#2}%
  \@ifundefined{subjclassname@#1}{%
    \ClassWarning{\@classname}{Unknown edition (#1) of Mathematics
      Subject Classification; using '1991'.}%
  }{%
    \@xp\let\@xp\subjclassname\csname subjclassname@#1\endcsname
  }%
}
\makeatother
\newcounter{constant}
\newcommand\ind{\hspace{5mm}}

\DeclareMathOperator\res{res}

\renewcommand\Re{\operatorname{Re}}
\renewcommand\Im{\operatorname{Im}}

	\renewcommand\L{\mathcal{L}}	\newcommand\D{\mathcal{D}}
\newcommand\CS{\mathrm{CS}}	\newcommand\Fib{\mathrm{Fib}}
\newcommand\M{\mathcal{M}}
\newcommand\oscp{{\bf p}}	
\newcommand\absc{{\sigma_c}}

\newcommand{\Flow}{\mathcal{F}}
\newcommand{\weight}{\mathfrak{w}}
\newcommand{\tweight}{\weight_{\text{tot}}}
		\newcommand{\orb}{\mathfrak{p}}

\newcommand\Int{\mathbb{Z}}	\newcommand\Nat{\mathbb{N}}
\newcommand\Real{\mathbb{R}}	\newcommand\Com{\mathbb{C}}
\newcommand\Fin{\mathbb{F}}
\newcommand\A{\mathbb{A}}

\newcommand\thick{\linethickness{1.1pt}}
\renewcommand\bar{\put(0,-5){\line(0,1){0}}}

\newtheorem{theorem}{Theorem}[section]

\theoremstyle{definition}
\newtheorem{definition}[theorem]{Definition}
\newtheorem{example}[theorem]{Example}

\newtheorem{conjecture}[theorem]{Conjecture}
\newtheorem{problem}[theorem]{Problem}

\theoremstyle{remark}
\newtheorem{remark}[theorem]{Remark}

\numberwithin{equation}{section}

\begin{document}

\title{Fractality, Self-Similarity and Complex Dimensions}
\author{Michel L.~Lapidus}
\address{Michel L.~Lapidus,
Department of Mathematics,
University of California, Riverside,
CA 92521-0135, USA}
\email{lapidus@math.ucr.edu}
\thanks{This work was partially supported by the US National Science Foundation under the research grant DMS-0070497 (for M.L.L.).}
\author{Machiel van Frankenhuijsen}
\address{Machiel van Frankenhuijsen,
Department of Mathematics,
Utah Valley State College,
Orem, UT 84058-5999, USA}
\email{vanframa@uvsc.edu}
\dedicatory{To Beno\^\i t Mandelbrot, on the occasion of his jubilee.}


\begin{abstract}
We present an overview of a theory of complex dimensions of self-similar fractal strings,
and compare this theory to the theory of varieties over a finite field from the geometric
and the dynamical point of view.
Then we combine the several strands to discuss a possible approach 
to establishing a cohomological interpretation of the complex dimensions.
\end{abstract}

\maketitle

\section{Introduction}

In this paper,
we survey some aspects of the theory of complex (fractal) dimensions of fractal strings,
as developed by the authors in the recent research monograph
{\em Fractal Geometry and Number Theory\/}~\cite{book} and pursued in the later
papers~[\citen{pot}--\ref{cxdssc}].
This theory builds upon
earlier work of the first author and his collaborators,
including Carl Pomerance,
Helmut Maier and Christina He
[\citen{mL91}--\ref{mL93c},\,\ref{mL95c},\,\citen{mLcP93}--\ref{mLcP96c},\,\citen
{LapMa},\,\citen{HeLap}].
In the present paper,
we stress the important case of self-similar strings
(which are one-dimensional self-similar geometries)
 along with the geometric and dynamical motivations and applications.
For the spectral theory,
in particular for a
 discussion of the question (\`a la Mark Kac~\cite{Kac})
{\em Can one hear the shape of a fractal drum?,\/}
going back to the
 Weyl--Berry Conjecture~\mbox{[\citen{hW12a},\,\citen{mB79}]},
we refer the interested reader to the above works 
and the relevant references therein,
including~[\mbox{\citen{BroCa,jKmL93}--\ref{jKmL01c}},\,\citen{mL94},\,\ref{mL97c}];
see also~\cite{aT03}.

For pedagogical reasons,
our goal in \S\S\ref{S: fs}--\ref{S: dyn} is to illustrate the mathematical theory of
 complex dimensions of self-similar strings by means of a few representative examples and results,
along with a discussion of some extensions and open problems.
We do not aim for the greatest generality or the most complete statements.
The interested reader may wish to consult~[\citen{book}--\ref{cxdssc}] 
(and the relevant references therein) for more details and generality.
In~\S\ref{S: shifts} and~\S\ref{S: coho},
as well as in~\S\ref{S: ss finite} below,
we explore and clarify the emerging analogies between self-similar strings (and self-similar geometries) on the one hand,
and varieties over finite fields and arithmetic geometries on the other hand.
\medskip

We give an overview of the contents of this article.
In~\S\ref{S: fs},
we define the notion of fractal string,
and introduce the geometric zeta function and (fractal) complex dimensions of a fractal string.
In~\S\ref{S: ss},
we focus on self-similar fractal strings.
The geometric zeta function of a fractal string $\L$ with lengths $l_1,l_2,\dots$ is defined
 in Equation~(\ref{E: zetaL}).
The poles of this function are called the {\em complex dimensions\/} of $\L$.
For a string that is self-similar under $N$ scalings by factors $r_1,\dots,r_N$ (without gaps,
see~\S\ref{S: ss}),
 the complex dimensions are the solutions to the equation $r_1^\omega+\dots+r_N^\omega=1$
(see~(\ref{E: zeta ss}) and~(\ref{E: eq}) below).
If $r_j=r^{k_j}$ for some positive number~$r$ and certain positive integers~$k_j$ ($j=1,\dots,N$),
this reduces to a polynomial equation,
and it is a transcendental equation when such integers do not exist.
We introduce the dichotomy of lattice vs.\ nonlattice self-similar strings in~\S\ref{S: l nl} to make this distinction.
This distinction is much like that between rational and irrational real numbers.

The complex dimensions of a lattice string lie periodically on finitely many vertical lines
 $\Re s=$ constant,
separated by a fixed distance $\oscp=2\pi/\log r^{-1}$,
the {\em oscillatory period\/} of the lattice string (see \S\ref{S: l nl}).
Therefore,
 we view a lattice string as an analogue of a variety over a finite field,
as explained below.
The space of lattice strings is dense in the space of all self-similar strings,
and nonlattice strings are shown to have properties that compel us to regard them as analogues of `infinite dimensional varieties'.
Moreover,
nonlattice strings with a very large number of scaling ratios form the bridge between
 self-similar and arithmetic geometries.
Some evidence (partly conjectural) for this is presented in~\S\ref{S: dfree}.
In that section,
we apply techniques from Diophantine approximation to obtain `dimension-free' regions for nonlattice strings.
The size of these dimension-free regions depends on how well
the ratios $\log r_j/\log r_1$ ($j=2,\dots,N$) can be approximated by rational numbers.
In our theory,
we deduce from the results in \S\ref{S: dfree} combined with our `explicit formulas'
 from~\cite{book} suitable asymptotic expansions with good error term for various geometric
 (\S\ref{S: l nl})
 or dynamical (\S\ref{S: dyn}) quantities associated with self-similar fractal strings:
the larger the dimension-free region,
the better the error term.

In~\S\ref{S: dyn},
we explain how the logarithmic derivative of the geometric zeta function of a fractal string
can be obtained as the generating function of the periodic
orbits of a continuous-time self-similar dynamical system.
We thus find an Euler-type product for the geometric zeta function.
This is analogous to the discrete-time dynamical system on a variety over a finite field,
namely,
the Frobenius flow,
as explained below.
Conjecturally~[\citen{Deninger},\,\citen{mystery}],
the corresponding arithmetic dynamical system is the shift on the real line (with infinitesimal generator the derivative).
In~\S\ref{S: shifts},
we propose an analysis of self-similar strings by measuring the overlap under shifts.
Finally,
in~\S\ref{S: coho} we present three tables,
 summarizing three aspects of these analogies.
\medskip

As was alluded to above,
whereas the work described in \S\S\ref{S: fs}--\ref{S: dyn} is completely rigorous and is part of a well-developed (and still growing) mathematical theory of complex dimensions of fractal strings,
 the content of \S\S\ref{S: shifts}--\ref{S: coho} is of a more speculative nature and should,
in the long-term,
be the object of further investigations and developments.
It builds upon well-known analogies in arithmetic (or Diophantine) geometry
(see,
e.g.,
\cite{aW41,Ka,aPiS95,Deninger}) and upon new ones proposed in~\mbox{[\citen{book},\,\S10.5]} and
 in the first author's forthcoming book~\cite{mL02},
{\em In Search of the Riemann Zeros}.

\subsection{Self-Similar and Finite Geometries}
\label{S: ss finite}

We close this introduction by providing some motivations and background material on varieties over
 a finite field.
We want to reassure the reader mainly interested in fractal geometry that no prior knowledge of arithmetic or algebraic geometry is required to read this paper.

The analogy between the geometric zeta function of self-similar fractal strings and the zeta function of a variety over a finite field
becomes apparent in the following simplest example.
The zeta function of the affine line\footnote
{Here,
and henceforth,
$p$ is a prime number and $\Fin_p$ denotes the finite field with $p$ elements.} over $\Fin_p$ is
defined as~$
\zeta_{\A^1/\Fin_p}(s)=\left(1-p\cdot p^{-s}\right)^{-1}$.
We compare this with the geometric zeta function of the Cantor string,
to be introduced in
Example~\ref{E: mCS},
\begin{gather}\label{E: zeta CS}
\zeta_\CS(s)=\frac1{1-2\cdot 3^{-s}}.
\end{gather}
Writing $D=\log_32=\log2/\log3$ for the Minkowski dimension of the Cantor string,
as will be explained in~\S\ref{S: fs},
we find
$
\zeta_{\A^1/\Fin_3}(s)=\zeta_\CS(s+D-1)$
or
$\zeta_{\A^1/\Fin_2}(s)=\zeta_\CS(sD)$.
These zeta functions have one line of poles,
at $\omega=1+ik\oscp$
($k\in\Int$ and $\oscp=2\pi/\log2$ or $\oscp=2\pi/\log3$,
respectively).

In general,
one defines the zeta function of a variety over a finite field as the generating function 
(or Mellin transform) of the counting measure of the positive divisors on the variety.
As such,
it is immediately clear that the zeta function can be obtained as an Euler product of factors that are defined in terms of the 
prime divisors of the variety.
A variety over a finite field comes equipped with an action of the Frobenius endomorphism.\footnote
{which is induced on the variety by the automorphism $x\mapsto x^p$ of $\Fin_p$.}
This defines a discrete-time flow on the variety,
the orbits of which are conjugacy classes of points on the variety,
which (for curves) are the prime divisors of the variety.
The logarithmic derivative of the Euler product is the generating function of the counting measure of the orbits of Frobenius.

One of the most important developments in the theory of algebraic varieties was the definition and subsequent development of a
cohomology theory.
Indeed,
the so-called \'etale cohomology,
which captures the combinatorics of families of \'etale covers of the variety,
provides a theory that can be compared to the classical singular homology and cohomology theories of manifolds.
Since the \'etale theory is defined purely algebraically,
it allows application to varieties defined over a finite field.
In particular,
one recovers the zeta function of the variety as the alternating product of the characteristic polynomials of the induced action 
of Frobenius on the \'etale cohomology groups.
The poles of this zeta function are located on `integer vertical lines' $\Re s=0,1,\dots,n$,
where $n$ is the dimension of the variety,
and the zeros are located on `half-integer vertical lines' $\Re s=\frac12,\frac32,\dots,n-\frac12$.
(See,
e.g.,~[\citen{Ka},\,\citen{FreKie},\,\citen{Deninger},\,\citen{aW41}--\ref{aW48c},\,\ref
{aW80c},\,\citen{aPiS95},\,I] for further information
about this beautiful subject.)

This theory was first modelled on the theory of the Riemann zeta function,
which is the first example of the zeta function of an `arithmetic geometry',
namely,
the spectrum of $\Int$.
There exist only one-dimensional\footnote
{Indeed,
the product of the spectrum of $\Int$ with itself has coordinate ring $\Int\otimes\Int=\Int$,
hence the product reduces to the diagonal~\cite{Haran}.} arithmetic geometries.
The (completed) Riemann zeta function has simple poles at $s=1$ and $s=0$,
hence one could say (by analogy with the case of a curve over a finite field) that the cohomology groups 
$H^0$ and $H^2$ are one-dimensional;
i.e.,
the spectrum of~$\Int$ (completed at the archimedean place) is connected and one-dimensional.
On the other hand,
if the Riemann hypothesis holds,
then the middle cohomology group $H^1$ could possibly be defined,
and should be infinite dimensional.
Moreover,
the logarithmic derivative of the Riemann zeta function is the generating function of the
 prime-power counting function.
It is not known how the prime numbers $p$ can be viewed as the (primitive) periodic orbits (of length~$\log p$) of a flow,
but there is some indication
that such a flow,
if it exists,
should take into account the smooth (archimedean) structure~\cite{Deninger,mystery}.
The simplest such flow is the shift on the real line.

\begin{remark}
In a book in preparation ([\citen{mL02}],
announced in~\cite{mL03}),
the first author has proposed a suitable notion of `quantized fractal string',
coined `fractal membrane',
and has significantly extended the above analogy between self-similar geometries and arithmetic geometries.
In particular,
lattice strings (or membranes) correspond to varieties over finite fields whereas nonlattice strings (or membranes) are viewed as
a counterpart of number fields.\footnote
{finite extensions of the field of rational numbers;
see,
e.g.,~[\citen{aPiS95},\,I\,\&\,II].}
Moreover,
the scaling ratios of self-similar geometries (or of `self-similar membranes'~\cite{mL02}) play the role of the `generalized primes'
attached to fractal membranes,
while expressions like~(\ref{E: zeta ss}) for the geometric zeta function of a self-similar string
(with possibly infinitely many scaling ratios and gaps)
correspond to the Euler product representation
(of the same nature as that for the classical Riemann zeta function~\cite{In,Pat,Ti})
of the partition function of a fractal membrane (as obtained in~\cite{mL02}).
Finally,
in joint work of the first author with Ryszard Nest~\cite{mLrN02},
it was recently shown that the fractal membranes (resp.\ self-similar membranes) introduced in~\cite{mL02} can be rigorously
constructed as the second quantization of fractal strings by using Fermi--Dirac
 (resp.,
Gibbs--Boltzmann or `free') statistics along with aspects of the theory of operator algebras and of 
Connes' noncommutative geometry~\cite{Connes}.
\end{remark}

\section{Fractal Strings}
\label{S: fs}

A {\em fractal string\/} is a one-dimensional drum with fractal boundary.
Thus,
a fractal string is a bounded open subset $\L$ of the real line.
We assume here that the connected components `touch'\footnote
{We might say,
the fractal string is {\em irreducible\/} or {\em connected,\/}
in the sense of irreducible schemes or connected manifolds.} each other;
i.e.,
the closure of $\L$ is a bounded connected interval $I$,
and the boundary $\partial\L$ of $\L$ is a disconnected fractal subset
of~$\Real$ with measure zero.
The lengths of the connected components (i.e.,
open intervals) of the string are called the {\em lengths\/} of $\L$
(a better terminology would be to speak of the {\em strings\/} of a {\em fractal harp;\/}
see~[\citen{book},\,p.\,8]).
We can list the lengths of a fractal string in nonincreasing order
$$
\L\colon l_1\geq l_2\geq l_3\geq\dots,
$$
each repeated according to its multiplicity.
The {\em geometric zeta function\/} of the fractal string is determined by the sequence $\L:$
\begin{gather}\label{E: zetaL}
\zeta_\L(s)=\sum_{j=1}^\infty l_j^s.
\end{gather}
Thus,
for example,
$\zeta_\L(1)=|I|$ is the length of the closure of the fractal string in the real line.
On the other hand,
$\zeta_\L(0)$ is the number of lengths of $\L$,
and we exclude the case when this number is finite.
Recall that the abscissa of convergence of a Dirichlet series is defined as the smallest real number $\absc$ such that the 
series~(\ref{E: zetaL}) converges (absolutely) in the half-plane $\Re s>\absc$.
Thus the abscissa of convergence of $\zeta_\L$ satisfies $0\leq \absc\leq1$.

We now compute the Minkowski (or box) dimension of the complement of the string in its interval.
(See also~[\citen{mL91}--\ref{mL92c},\,\citen{mLcP93}].)
When the two endpoints of an interval of length $l_j$ are covered by intervals (or disks) of radius $\varepsilon$,
then these disks overlap if $l_j<2\varepsilon$,
covering a length $l_j$,
or they do not overlap if $l_j\geq2\varepsilon$,
in which case they cover a length of $2\varepsilon$.
Moreover,
there are the two pieces of length~$\varepsilon$ sticking out to the left and the right of $I$.
Thus the length $V(\varepsilon)=V_\L(\varepsilon)$ of the (two-sided) $\varepsilon$-neighborhood of the string is given by\footnote
{For two subsets of $\Real$,
 $A+B$ denotes the set $\{a+b\colon a\in A,b\in B\}$.
Thus,
$A+(-\varepsilon,\varepsilon)$ is the open $\varepsilon$-neighborhood of $A$.
Also,
$|A|$ denotes Lebesgue measure in $\Real$.\label{F: set sum}}
\begin{gather}\label{E: Ve}
V(\varepsilon)=|\partial\L+(-\varepsilon,\varepsilon)|
=2\varepsilon+\sum_{j\colon l_j\geq2\varepsilon}2\varepsilon+\sum_{j\colon l_j<2\varepsilon}l_j.
\end{gather}

\begin{remark}\label{R: two-sided}
In~[\citen{arith}--\ref{cxdssc}],
we used the inner neighborhood of the string,
to allow for strings with lengths not necessarily embedded in the real line.
The present approach is less intrinsic,
but more suitable to study properties of self-similarity
(also see the end of~\S\ref{S: cxd}).
\end{remark}

The {\em Minkowski dimension $D$\/} of $\partial\L$ is the infimum
of the numbers $d$ such that~$V(\varepsilon)=o(\varepsilon^{1-d})$ as $\varepsilon\to0^+$:
$
D=\inf\left\{d\geq0\colon\varepsilon^{d-1}V(\varepsilon)\to0\text{ as }\varepsilon\to0^+\right\}
$.
Further,~$\partial\L$
 is {\em Minkowski measurable\/} if the following limit exists in $(0,\infty)$,
\begin{gather}\label{E: M}
\M=\M(D;\L):=\lim_{\varepsilon\to0^+}\varepsilon^{D-1}V(\varepsilon).
\end{gather}
(See Remark~\ref{R: ss boundary} below and,
e.g.,
[\citen{mL91}--\ref{mL93c},\,\citen{mLcP93}--\ref{mLcP96c},\,\citen{kF90},\,Ch.\,3,\,\citen
{Mat},\,Ch.\,5].)

We now show that the Minkowski dimension can be recovered as the abscissa of convergence of
 $\zeta_\L$.
(This was first observed by the first author in~[\citen{mL92}--\ref{mL93c}] using a result of Besicovich and Taylor~\cite{BT};
see also~[\citen{mLcP93},\,\citen{LapMa}].)

\begin{theorem}
The abscissa of convergence $\absc$ of the geometric zeta function\/~$\zeta_\L$ coincides with $D$,
the Minkowski dimension of\/ $\partial\L$.
\end{theorem}

\begin{proof}
Let $s>d>D$.
There exists a constant $\refstepcounter{constant}\label{C: vol}C_{\ref{C: vol}}$ 
such that $V(\varepsilon)\leq C_{\ref{C: vol}}\varepsilon^{1-d}$.
For $\varepsilon=l_n/2$,
we obtain
$
(n+1)l_n\leq (n+1)l_n+\sum_{j=n+1}^\infty l_j=V(l_n/2)\leq C_{\ref{C: vol}}(l_n/2)^{1-d}.
$
It follows that 
$\refstepcounter{constant}\label{C: ln}
l_n^s\leq C_{\ref{C: ln}} (n+1)^{-s/d},
$
for some constant $C_{\ref{C: ln}}$.
Hence the series~(\ref{E: zetaL}) converges for $s>d$,
so that $\absc\leq d$.
Since this holds for every~$d>D$,
we obtain~\mbox{$\absc\leq D$}.
If $\absc=1$,
we conclude that $D=\absc$,
since $V(\varepsilon)\leq2\varepsilon+|I|$ is bounded.

Otherwise,
let $\absc<s<1$.
Then the series~(\ref{E: zetaL}) converges.
Since the sequence of lengths is nonincreasing,
we find that
$\refstepcounter{constant}\label{C: sum}
nl_n^s\leq\sum_{j=1}^nl_j^s\leq C_{\ref{C: sum}}.
$
Hence~$l_n\leq (C_{\ref{C: sum}}/n)^{1/s}$.
Given~\mbox{$\varepsilon>0$},
it follows that $l_n<2\varepsilon$ for~$n>C_{\ref{C: sum}}(2\varepsilon)^{-s}$.
For $j\leq C_{\ref{C: sum}}(2\varepsilon)^{-s}$,
we estimate the $j$-th term (in the first or the second sum) in~(\ref{E: Ve}) by $2\varepsilon$,
and for~\mbox{$j>C_{\ref{C: sum}}(2\varepsilon)^{-s}$},
 we estimate this term in the second sum by $(C_{\ref{C: sum}}/j)^{1/s}$.
Thus we find
$\refstepcounter{constant}V(\varepsilon)\leq C_\theconstant(2\varepsilon)^{1-s}$,
so that~$D\leq s$.
Since this holds for every real number~$s>\absc$,
we conclude that~$D\leq\absc$.
Hence these two quantities coincide.
\end{proof}

\subsection{Complex Dimensions of a Fractal String}
\label{S: cxd}

We assume that $\zeta_\L$ has a meromorphic continuation to a neighborhood of the half-plane
 $\Re s\geq D$.
Clearly then,
$D$ is a pole of $\zeta_\L$.
The geometric zeta function may have other singularities,
which are then located on or to the left of the line $\Re s=D$.
These singularities are necessarily poles,
and they are called the {\em complex dimensions\/} of $\L$.
We denote by $\D=\D_\L$ the set of (visible) complex dimensions of $\L$ (within the given region).
Clearly,
it is a discrete subset of $\Com$,
and therefore is at most countable.

The complex dimensions enter into the explicit formulas for the various geometric,
dynamical,
 and spectral quantities related to $\L$.
These (pointwise or distributional) explicit formulas---in the sense of number
 theory~\cite{In,Pat,Ti,aW52},
but more general---are established in~[\citen{book},\,Ch.\,4] and then applied and adapted to a number of situations throughout the rest of that book,
especially in~[\citen{book},\,Chs.\,5--7,\,9]
(see also~[\citen{pot}--\ref{cxdssc}],
along with the forthcoming second edition of~\cite{book}).

As an illustration,
 we have the explicit formula for the volume (i.e.,
length) of the tubular neighborhood~$V(\varepsilon)$,
which we formulate here in the case when the complex dimensions are simple.
See~[\citen{book},\,Ch.\,6] for the general formulation and additional details.

\begin{theorem}
Let\/ $\zeta_\L(s)$ have a meromorphic continuation to a neighborhood of\/
 $\Re s\geq\sigma_0$ (for some $\sigma_0\leq D$),
with simple poles and such that\/ $s=0$ is not a pole of\/ $\zeta_\L$.
Then,
under mild growth conditions on $\zeta_\L,$
\begin{gather}\label{E: Ve exp}
V(\varepsilon)=\sum_\omega\res\left(\zeta_\L;\omega\right)
\frac{(2\varepsilon)^{1-\omega}}{\omega(1-\omega)}+2\varepsilon(1+\zeta_\L(0))+O\left(\varepsilon^{\sigma_0}\right)
\text{ as }\varepsilon\to0^+.
\end{gather}
The sum extends over all poles of\/ $\zeta_\L$ in $\sigma_0\leq\Re\omega\leq D$
(i.e.,
over all `visible' complex dimensions $\omega$ of\/ $\L$).
\end{theorem}

This formula expresses $V(\varepsilon)$ as a sum of oscillatory terms of the form $C_\omega\varepsilon^{1-\omega}$,
where~$\omega$ runs over all (visible) complex dimensions of $\L$.
These terms grow asymptotically as $\varepsilon^{\Re\omega}$ in size,
and have a multiplicative period $r=\exp(2\pi/\Im\omega)$;
i.e.,~$\varepsilon^{\omega}$ and~$(r\varepsilon)^{\omega}$ have the same argument.
The coefficient of each oscillatory term depends in particular on $\res(\zeta_\L;\omega)$,
the residue of the geometric zeta function at~$\omega$.
By analogy with Weyl's formula for the volume of the $\varepsilon$-neighborhoods of a (smooth) submanifold of Euclidean space
([\citen{hW39},\,\citen{mBbG88},\,\S6.6--6.9,\,esp.\,Thm.\,6.9.9]),
we expect eventually to be able to obtain a suitable geometric interpretation of these values;
 see~[\citen{book},\,\S6.1.1].
For example,
in view of~(\ref{E: M}),
the residue at~$D$ is closely related to the Minkowski content (see~[\citen{book},\,Thm.\,6.12]):
\begin{gather}\label{E: MDL}
\M(D;\L)=\res(\zeta_\L;D)\frac{2^{1-D}}{D(1-D)},
\end{gather}
provided that there are no other complex dimensions with real part $D$.
Moreover,
this is an important reason why we prefer to work here with two-sided (rather than with one-sided) neighborhoods of $\L$
(see Remark~\ref{R: two-sided} above).
Indeed,
in the former case,
the geometric interpretation of Weyl's tube formula~\cite{hW39} is valid for submanifolds of any dimension,
whereas in the latter case,
it is only valid for even-dimensional submanifolds;
see,
e.g.,
\cite{mBbG88},
{\em loc.\,cit.}

\section{Self-Similar Fractal Strings}
\label{S: ss}

\begin{figure}[tb]
\begin{picture}(324,105)(-15,10)
\unitlength 0.4pt
\put(0,260){\put(364.5,20){\makebox(0,0){$L=3$}}
\bar{\thick\line(1,0){729}}\bar
}

\put(0,200){
\put(0,0){\line(1,0){729}}
\put(40.5,20){\makebox(0,0){$r_1=\frac19$}}\bar{\thick\line(1,0){81}}
\bar\put(81,0){\put(40.5,20){\makebox(0,0){$r_2=\frac19$}}\bar{\thick\line(1,0){81}}\bar
	\put(243,0){\put(40.5,20){\makebox(0,0){$r_3=\frac19$}}
\bar{\thick\line(1,0){81}}\bar
\put(81,0){
\put(13.5,20){\makebox(0,0){$r_4$}}\bar{\thick\line(1,0){27}}\bar
		\put(27,0){\put(13.5,20){\makebox(0,0){$r_5$}}\bar{\thick\line(1,0){27.3}}\bar}}}}
}

\put(0,140){
\put(0,0){\line(1,0){81}}
\bar{\thick\line(1,0){9}}\bar\put(9,0){\bar{\thick\line(1,0){9}}\bar
	\put(27,0){\bar{\thick\line(1,0){9}}\bar\put(9,0){\bar{\thick\line(1,0){3}}\bar
		\put(3,0){\bar{\thick\line(1,0){3}}\bar}}}}
\put(153,0){
\put(0,0){\line(1,0){81}}
\bar{\thick\line(1,0){9}}\bar\put(9,0){\bar{\thick\line(1,0){9}}\bar
	\put(27,0){\bar{\thick\line(1,0){9}}\bar\put(9,0){\bar{\thick\line(1,0){3}}\bar
		\put(3,0){\bar{\thick\line(1,0){3}}\bar}}}} }
\put(477,0){
\put(0,0){\line(1,0){81}}
\bar{\thick\line(1,0){9}}\bar\put(9,0){\bar{\thick\line(1,0){9}}\bar
	\put(27,0){\bar{\thick\line(1,0){9}}\bar\put(9,0){\bar{\thick\line(1,0){3}}\bar
		\put(3,0){\bar{\thick\line(1,0){3}}\bar}}}} }
\put(639,0){
\put(0,0){\line(1,0){27}}
\bar{\thick\line(1,0){3}}\bar\put(3,0){\bar{\thick\line(1,0){3}}\bar
	\put(9,0){\bar{\thick\line(1,0){3}}\bar\put(3,0){\bar{\thick\line(1,0){1}}\bar
		\put(1,0){\bar{\thick\line(1,0){1}}\bar}}}} }
\put(693,0){
\put(0,0){\line(1,0){27}}
\bar{\thick\line(1,0){3}}\bar\put(3,0){\bar{\thick\line(1,0){3}}\bar
	\put(9,0){\bar{\thick\line(1,0){3}}\bar\put(3,0){\bar{\thick\line(1,0){1}}\bar
		\put(1,0){\bar{\thick\line(1,0){1}}\bar}}}} }
}

\put(124,100){\makebox(0,0){$\vdots$}}\put(526.5,100){\makebox(0,0){$\vdots$}}\put(688.5,100){\makebox(0,0){$\vdots$}}

\put(0,55){
\put(81,0){\put(40.5,-20){\makebox(0,0){$g_1=\frac19$}}\line(1,0){81}
	\put(81,0){\put(121.5,-20){\makebox(0,0){$g_2=\frac13$}}\put(121.5,20){\makebox(0,0){$\L$}}
	\line(1,0){243}\put(81,0){\put(40.5,-20){\makebox(0,0){$g_3=\frac19$}}
	\line(1,0){81}\put(27,0){\put(13.5,-20){\makebox(0,0){$g_4$}}\line(1,0){27}}}}}
\put(0,0){
\put(9,0){\line(1,0){9}\put(9,0){\line(1,0){27}\put(9,0){\line(1,0){9}\put(3,0){\line(1,0){3}}}}}}
\put(162,0){
\put(9,0){\line(1,0){9}\put(9,0){\line(1,0){27}\put(9,0){\line(1,0){9}\put(3,0){\line(1,0){3}}}}}}
\put(486,0){
\put(9,0){\line(1,0){9}\put(9,0){\line(1,0){27}\put(9,0){\line(1,0){9}\put(3,0){\line(1,0){3}}}}}}
\put(648,0){
\put(3,0){\line(1,0){3}\put(3,0){\line(1,0){9}\put(3,0){\line(1,0){3}\put(1,0){\line(1,0){1}}}}}}
\put(702,0){
\put(3,0){\line(1,0){3}\put(3,0){\line(1,0){9}\put(3,0){\line(1,0){3}\put(1,0){\line(1,0){1}}}}}}}
\end{picture}
\caption{The modified Cantor string,
with five scaling ratios $r_1=r_2=r_3=\frac{1}{9}$, $r_4=r_5=\frac{1}{27}$,
and four gaps $g_1=g_3=\frac19$, $g_2=\frac13$, $g_4=\frac1{27}$.}
\label{F: ss string}
\end{figure}
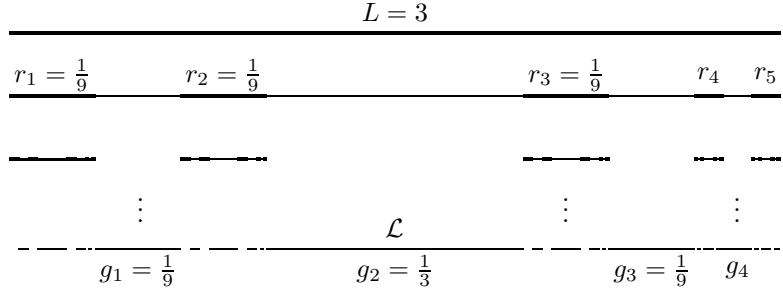

Let positive {\em scaling ratios\/} $r_1\geq r_2\geq\dots\geq r_N$ and
 {\em gaps\/} $g_1,\dots,g_K$ be given
 such that~$r_1+\dots+r_N+g_1+\dots+g_K=1$.
Figure~\ref{F: ss string} explains how one constructs a self-similar fractal string with an `initiator' $r_1,\dots,r_N,g_1,\dots,g_K$
in an initial interval of length $L$.
(See~[\citen{mF03},\,\citen{cxdss}].)
The initial interval is replaced by scaled copies,
leaving gaps in between.
The gaps and scaled gaps form the lengths of the string $\L$,
which are the products of the form $Lg_kr_{j_1}\dots r_{j_q}$,
$1\leq k\leq K$,
$1\leq j_1,\dots,j_q\leq N$,
$q\geq0$.
The multiplicity of a length $l$ is the number of such products with value~$l$.
In~[\citen{book},\,\ref{cxdssc}],
we have computed the geometric zeta function of $\L$,
\begin{gather}\label{E: zeta ss}
\zeta_\L(s)=\frac{L^s\sum_{k=1}^Kg_k^s}{1-\sum_{j=1}^Nr_j^s}.
\end{gather}
It follows that the geometric zeta function of a self-similar string has a meromorphic continuation to the entire complex plane,
given by~(\ref{E: zeta ss}).
Further,
the set $\D_\L$ of complex dimensions of $\L$ is a subset of the set of solutions of the equation
\begin{gather}\label{E: eq}
r_1^\omega+r_2^\omega+\dots+r_N^\omega=1\qquad\omega\in\Com.
\end{gather}
It may be a proper subset if zeros of the denominator are cancelled by zeros of the numerator of $\zeta_\L$,
as we show in this paper
(see Examples~\ref{E: mCS} and~\ref{E: mFib} below).
Moreover,
$D$ is the unique {\em real\/} solution to~(\ref{E: eq}).
Note that $D$ is always a simple root of this equation,
and that it is never cancelled by a root of the numerator of~$\zeta_\L$,
since the numerator of (\ref{E: zeta ss}) does not have real roots.

Of particular interest is the case where there is only one gap,
$K=1$,
necessarily of size
$g=1-r_1-\dots-r_N$.
Then,
$\zeta_\L(s)=(gL)^s/(1-\sum_{j=1}^Nr_j^s)$,
and the complex dimensions are exactly the (complex) solutions of Equation~(\ref{E: eq}).
This case was studied in~\cite{book};
see also~\S\ref{S: dyn}.
Another natural example is when the scaled copies of $I$ are interspersed by gaps,
so that
$K=N-1$.
Only in this case,~
\mbox{$\zeta_\L(0)=-1$},
so that there is no term of order $\varepsilon$ in~(\ref{E: Ve exp}).
In general,~\mbox{$\zeta_\L(0)=K/(1-N)$}.

\begin{remark}\label{R: ss boundary}
After a possible rearrangement,
any self-similar set in $\Real$ (without too much overlap\footnote
{that is,
satisfying the so-called `open set condition' [\citen{jH81},\,\citen{kF90},\,p.\,112].})
can be viewed as the boundary $\partial\L$ of a self-similar string~$\L$ (possibly with multiple gaps);
and conversely.
Further,
it is well-known that for such self-similar sets,
the Minkowski and the Hausdorff dimensions coincide.
Moreover,
since $D$ is the unique real solution of the
 `complexified Moran equation'~(\ref{E: eq})~\cite{Mor},
it also coincides with the similarity dimension of $\partial\L$
 (in the sense of Mandelbrot~\cite{Mandelbrot});
see,
e.g.,~[\citen{kF90},\,Thms.\,9.1,\,9.3].
\end{remark}

We next give the example of the simplest self-similar string.
The second example shows that the initiator of a self-similar string is in general not unique.

\begin{example}[Cantor and Modified Cantor String]\label{E: mCS}
The string with five scaling ratios $r_1=r_2=r_3=1/9$,
$r_4=r_5=1/27$ and four gaps $g_1=1/3$,
$g_2=g_3=1/9$,
$g_4=1/27$ in an interval of length $3$ (see Figure~\ref{F: ss string}),
has for geometric zeta function
\begin{gather}\label{E: zeta MCS}
\zeta_\L(s)=3^s\frac{3^{-s}+2\cdot3^{-2s}+3^{-3s}}{1-3\cdot3^{-2s}-2\cdot3^{-3s}}.
\end{gather}
The denominator factors as follows:
$1-3x^2-2x^3=(1-2x)(1+x)^2$,
and the numerator as $x(1+x)^2$,
with $x=3^{-s}$.
Hence the zeta function simplifies to the geometric zeta function~(\ref{E: zeta CS})
 of the Cantor string,
 with two scaling ratios $r_1=r_2=1/3$ and one gap $g_1=1/3$,
in an interval of length $3$.
Thus the sequence of its lengths coincides with the sequence of lengths of the connected
 components of the complement of the Cantor set,
see~[\citen{book},\,Ex.\,2.2.1].
The poles (of~(\ref{E: zeta MCS}) or~(\ref{E: zeta CS})) are simple,
 located at $D+ik\oscp$ ($k\in\Int$),
with $D=\log_32$ and $\oscp=2\pi/\log 3$;
the residue at each pole is equal to $1/\log 3$.
In particular,
the complex dimensions of the Cantor string lie periodically on a single vertical line,
defined by $\Re s=D$.
\end{example}

\begin{remark}
The initiator of a self-similar fractal string is not unique 
 since one can always find a Dirichlet polynomial\/ $L^s(g_1^s+\dots+g_K^s)$ with positive
 coefficients such that the product of this Dirichlet polynomial with the denominator
 of\/ $\zeta_\L$ is again of the form $1-r_1^s-\dots-r_N^s$.
Indeed,
repeated application of the identity~$
\left(1-r_1^s-\dots-r_N^s\right)^{-1}=\left(1+r_1^s+\dots+r_N^s\right)/
\left(1-(r_1^s+\dots+r_N^s)^2\right)
$ 
shows that every fractal string has infinitely many initiators.
Example~\ref{E: mFib} below gives an alternative initiator for the Fibonacci string.
\end{remark}

The following example shows that
the complex dimensions of a self-similar fractal string are not necessarily simple
 (see~[\citen{book},\,Ex.\,2.2.3]).

\begin{example}[Multiple Complex Dimensions]\label{E: not simple}
The string with the same scaling ratios as in Example~\ref{E: mCS} but one gap $g_1=1/3+2/9+1/27$ has complex dimensions on a line
$\pi i/\log3+2k\pi i/\log3$ ($k\in\Int$),
each with multiplicity two.
\end{example}

\begin{example}[Fibonacci String]\label{E: Fib}
The Fibonacci string,
introduced in~[\citen{book},\,Ex.\,2.2.2],
 has the initiator $r_1=1/2$,
$r_2=g_1=1/4$,
in an interval of length~$4$.
Thus
the geometric zeta function of this string is
$\zeta_\Fib(s)=\left({1-2^{-s}-2^{-2s}}\right)^{-1}$.
This is a string with lengths $2^{-n}$,
of multiplicity $f_{n+1}$,
the $(n+1)$-st Fibonacci number (defined by $f_0=0$,
$f_1=1$ and $f_{n+1}=f_n+f_{n-1}$).
We find two (discrete) lines of complex dimensions:
$\{D+in\oscp\}_{n\in\Int}$,
with $D=\log_2\phi$ and $\oscp=2\pi/\log2$,
and the line $\{-D-i\oscp/2+in\oscp\}_{n\in\Int}$.

The volume $V_\Fib(\varepsilon)$
of the tubular neighborhood of the Fibonacci string can be computed directly,
thereby illustrating formula~(\ref{E: Ve}).
Let $\phi=(1+\sqrt{5})/2$ be the golden ratio and recall that 
$f_n=\left(\phi^n-(1-\phi)^n\right)/\sqrt{5}$.
By~(\ref{E: Ve}),
we have
$V_\Fib(\varepsilon)=
2\varepsilon\sum_{2^{-n}\geq2\varepsilon}f_{n+1}+\sum_{2^{-n}<2\varepsilon}f_{n+1}2^{-n}$.
By the above formula for $f_n$,
both sums are geometric.
We write $x=\log_2(2\varepsilon)^{-1}$,
so that $x$ increases by one unit if $\varepsilon$ is halved in value.
Evaluating the above sums,
we find that
\begin{align*}
V_\Fib(\varepsilon)&=\frac{(2\varepsilon)^{1-D}}{\sqrt5}
\left(\phi^3\phi^{-\{x\}}+\phi^4(\phi/2)^{-\{x\}}\right)\\&\hspace{5mm}
+\frac{(2\varepsilon)^{1+D}}{\sqrt5}(-1)^{[x]}
\left(\phi^{-3}\phi^{\{x\}}-\phi^{-4}(2\phi)^{\{x\}}\right).
\end{align*}
Here,
$[x]$ denotes the integer part,
and $\{x\}=x-[x]$ the fractional part of the real number~$x=\log_2(2\varepsilon)^{-1}$.
The two functions in parentheses are periodic (and continuous).
Computing their Fourier series,
we recover for the Fibonacci string the explicit formula~(\ref{E: Ve exp}),
even without error term:\footnote
{The general theory of~[\citen{book},\,Ch.\,4] also enables us to obtain this explicit formula
 without error term.
An analogous statement holds for all lattice strings
(see~\S\ref{S: l} and~[\citen{book},\,Thms.\,4.8,\,6.21]).}
\begin{equation}
\begin{aligned}
V_\Fib(\varepsilon)
&=\frac\phi{\sqrt5\log2}\sum_{k=-\infty}^\infty
\frac{(2\varepsilon)^{1-D-ik\oscp}}{(D+ik\oscp)(1-D-ik\oscp)}\\
&\hspace{5mm}
+\frac{\phi-1}{\sqrt5\log2}
\sum_{k=-\infty}^\infty
\frac{(2\varepsilon)^{1+D-i\oscp/2-ik\oscp}}{(-D+i\oscp/2+ik\oscp)(1+D-i\oscp/2-ik\oscp)}.
\end{aligned}
\end{equation}
\end{example}

\begin{example}[Modified Fibonacci String]\label{E: mFib}
The initiator $r_1=r_2=g_1=1/4$,
$r_3=g_2=1/8$ in an interval of length~$4$ generates a string with geometric zeta function
$
\zeta_\L(s)=2^{2s}\left(2^{-2s}+2^{-3s}\right)/\left(1-2\cdot2^{-2s}-2^{-3s}\right)
=\left(1-2^{-s}-2^{-2s}\right)^{-1}
$.
Hence the sequence of lengths of this string coincides with that of the Fibonacci string of Example~\ref{E: Fib}.
\end{example}

\section{The Lattice and the Nonlattice Case: Periodicity vs.\ Quasiperiodicity}
\label{S: l nl}

The lattice/nonlattice dichotomy comes from renewal
 theory~[\citen{wF66},\,Ch.\,XI].
It was used in a related context in~[\citen{sL88}--\ref{sL89c},\,\citen{rS90}--\ref{rS909393c}]
and~\cite{mL93};
see also,
e.g.,
[\citen{jKmL93,mL94,kF95}--\ref{kF97c},\,\citen{book,bHmL99}],
 and the relevant references therein.

\subsection{The Lattice Case}
\label{S: l}

The self-similar fractal string $\L$ is a {\em lattice\/} string if the multiplicative subgroup
$G=\prod_{j=1}^Nr_j^\Int$
 of $\Real^*_+$ generated by~$r_1,\dots,r_N$
has rank~$1$ as a free abelian group.
In that case,
 $G$ is a discrete subgroup $r^\Int$ of $\Real^*_+$.
The number $r\in(0,1)$ is called the {\em multiplicative generator of $\L$\/}.

Examples~\ref{E: mCS} and~\ref{E: not simple}--\ref{E: mFib} are all examples of lattice strings.
As we saw in the example of the Fibonacci string,
in the lattice case we can obtain complete information about the geometry of the string.
This is reflected in the following theorem,
which describes the structure of the complex dimensions of a lattice string
 (see~[\citen{book},\,Thm.\,2.13]).

\begin{theorem}[Lattice Case: Periodic Patterns]\label{T: l}
Let $\L$ be a lattice self-similar string with multiplicative generator $r$.
Then $\zeta_\L$ is a rational function of\/ $r^s,$
and hence there exist finitely many solutions $\omega_1:=D, \omega_2,\dots,\omega_m$ 
of Equation~(\ref{E: eq}) such that the complex dimensions of\/ $\L$ lie on finitely many vertical lines
\begin{gather}\label{E: line}
\omega=\omega_u+ik\oscp,\qquad k\in\Int,\quad u=1,\dots,m,
\end{gather}
where $\oscp=2\pi/\log r^{-1}$ is the {\em oscillatory period\/} of $\L$.
In particular,
the complex dimensions of a lattice string exhibit a\/ {\em periodic\/} pattern.

Since there are complex dimensions other than $D$ with real part equal to\/ $D$,
it follows from Equation~(\ref{E: Ve exp}) that\/ $\varepsilon^{D-1}V(\varepsilon)$ oscillates without a limit as $\varepsilon\to0^+$.
Hence a lattice string is not Minkowski measurable.
\end{theorem}

To prove this,
observe that in the lattice case,
there exist integers $1\leq k_1\leq k_2\leq\dots\leq k_N$ such that $r_j=r^{k_j}$
 for $j=1,\dots,N$,
so that Equation~(\ref{E: eq}) is equivalent to
\begin{gather}\label{E: lattice eq}
z^{k_1}+z^{k_2}+\dots+z^{k_N}=1,\qquad z=r^\omega.
\end{gather}
Hence $m$,
the number of lines of complex dimensions,
satisfies $m\leq k_N$.
Moreover,
some complex dimensions,
and even entire lines of complex dimensions,
can be canceled by roots of the numerator of~(\ref{E: zeta ss}).
See also Remark~\ref{R: nl gaps}.
\medskip

For example,
if $L=1$ and $\L$ has one gap,
then we obtain,
for $0<\varepsilon<(2r)^{-1}$,
\begin{equation}\label{E: vol tub lattice}
\begin{aligned}
V(\varepsilon )&=\res\left(\zeta_\L (s);D\right)
\sum_{n\in\Int}\frac{(2\varepsilon )^{1-D-in\oscp}}{(D+in\oscp )(1-D-in\oscp )}\\
&\hspace{20mm}
+\sum_{\Re\omega <D}\res\left(\frac{\zeta_\L (s)(2\varepsilon )^{1-s}}{s(1-s)};\omega\right)
+2\varepsilon\left(\zeta_\L (0)+1\right)\\
&=\res\left(\zeta_\L (s);D\right)(2\varepsilon )^{1-D}G\left(\log_r (2\varepsilon )\right)
+o\left(\varepsilon^{1-\Theta}\right) ,\mbox{ as }\varepsilon\to 0^+,
\end{aligned}
\end{equation}
where the line $\Re s=\Theta <D$ is the first line of complex dimensions that lies (strictly) to the left of $D$,
and $G$ is the periodic function
\begin{equation}\label{E: TM G}
\begin{aligned}
G(x)&=\sum_{n\in\Int}\frac{e^{2\pi inx}}{(D+in\oscp )(1-D-in\oscp )}
=\log r^{-1}\left( \frac{r^{D\{ x\}}}{1-r^D}+\frac{r^{(D-1)\{ x\}}}{r^{D-1}-1}\right).
\end{aligned}
\end{equation}
Since this periodic function is nonconstant,
it follows that $\L$ is not Minkowski measurable.

In general,
each (discrete) line of complex dimensions $\{\omega_u+ik\oscp\colon k\in\Int\}$
 (for~$u=1,\dots,m$) of a lattice string gives rise to a
(multiplicatively) periodic function times a suitable power of $\varepsilon$ (if $\omega_u$ is a simple pole\footnote
{In general,
times $\varepsilon^{1-\omega_u}P(\log\varepsilon^{-1})$,
where $P$ is a polynomial of degree one less than the multiplicity of $\omega_u$ as a pole of $\zeta_\L$.
Recall that $D$ is always a simple pole.})
in the corresponding explicit formula
(see~[\citen{book},\,\S6.3.1,\,esp.\,Thm.\,6.21,\,Ex.\,6.25]).
\medskip

\begin{figure}[p]
\begin{picture}(360,543)(0,5)		
\put(0,395){\put(60,70){\makebox(0,0){$\sqrt{2}\approx3/2$}}
	\put(66.,102){$\scriptstyle\oscp$}
	\epsfxsize 4.1cm\epsfbox{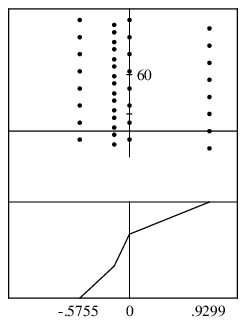}}
\put(120,395){\put(60,70){\makebox(0,0){$\sqrt{2}\approx7/5$}}
	\put(66.,115){$\scriptstyle\oscp$}
	\epsfxsize 4.1cm\epsfbox{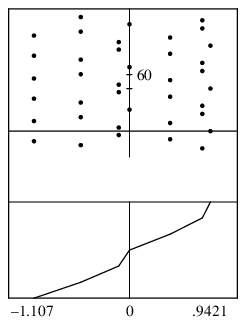}}
\put(240,395){\put(60,70){\makebox(0,0){$\sqrt{2}\approx17/12$}}
	\put(66.,146){$\scriptstyle\oscp$}
	\epsfxsize 4.1cm\epsfbox{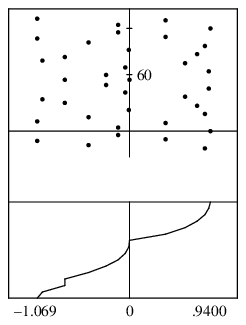}}
\put(0,0){\put(60,70){\makebox(0,0){$\sqrt{2}\approx41/29$}}
	\put(66,219.5){$\scriptstyle\oscp$}
	\epsfxsize 4.1cm\epsfbox{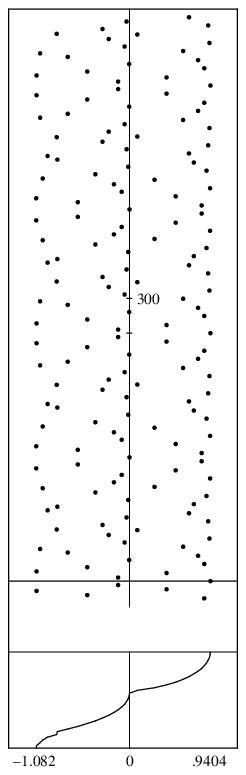}}
\put(120,0){\put(60,70){\makebox(0,0){$\sqrt{2}\approx99/70$}}
	\put(66,245.5){$\scriptstyle\oscp/2$}
	\epsfxsize 4.1cm\epsfbox{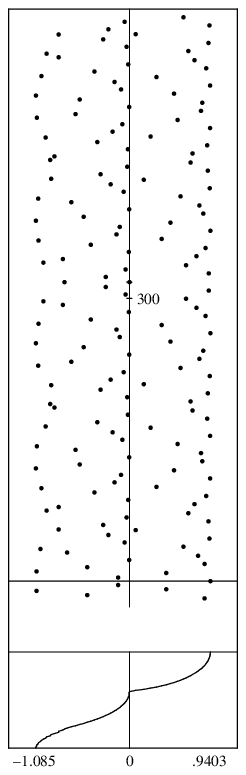}}
\put(240,0){\put(60,70){\makebox(0,0){$\sqrt{2}\approx239/169$}}
	\epsfxsize 4.1cm\epsfbox{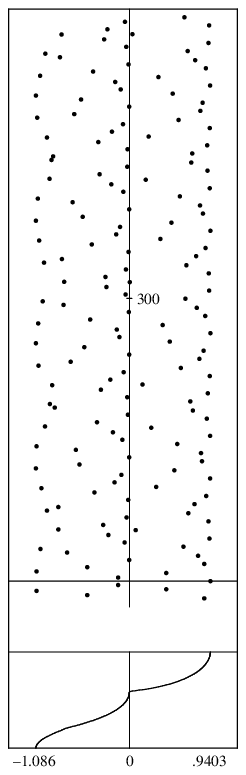}}
\end{picture}
\caption{Six stages of approximation of the complex dimensions of the nongeneric nonlattice string of Example~\ref{E: nong ex},
with $r_1=2^{-1}$,
$r_2=2^{-2}$ and $r_3=2^{-1-\sqrt{2}}$.
{\em Lower parts:} The density graph of the real parts of the complex dimensions.}
\label{F: stages}
\end{figure}

The oscillatory period $\oscp$ of a lattice string $\L$ can be interpreted as the generator of the 
{\em spectrum of self-similarity\/} of $\L$.
Indeed,
the lengths of $\L$ are of the form $r^ng_kL$ ($k=1,\dots,K$,
$n\in\Nat$),
with a multiplicity that grows like $r^{Dn}$;
i.e.,
exponentially with $n$.
Figure~\ref{F: stages} gives diagrams of the complex dimensions of six lattice strings that
 approximate the nonlattice string of Example~\ref{E: nong ex}.
For example,
the diagram labeled $\sqrt2\approx41/29$ gives the complex dimensions of the lattice string with
$r_1=r^{29}$,
$r_2=r^{2\cdot29}$,
and $r_3=r^{41+29}$,
where $r=2^{-1/29}$ (and one gap $g=1/4-r_3$).
Its oscillatory period is $\oscp=58\pi/\log 2\approx262.9$.

The lower parts of the diagrams in this figure give the relative density of the real parts of the
 complex dimensions within one period.
These density graphs were obtained in the following way.
For each solution $z$ of Equation~(\ref{E: lattice eq}),
we computed the real part of the complex dimension $\omega$,
namely
$\Re\omega=\log_r|z|$.
Then we ordered this set to obtain an increasing sequence $x_1\leq x_2\leq\dots\leq x_{k_N}$ of
 real parts,
and we plotted the points~$(x_j,j)$.
Thus,
steep parts of the density graph indicate that there are many complex dimensions with almost the same real part,
whereas flat parts indicate that the sequence $x_j,x_{j+1},\dots$ increases relatively fast.

\subsection{The Nonlattice Case}
\label{S: nl}

The {\em nonlattice case\/} is when the group $G$ defined at the beginning of~\S\ref{S: l} is
 dense in~$\Real^*_+$.
Moreover,
the fractal string $\L$ is said to be a {\em generic nonlattice\/} string if the cardinality 
$M=\#\{r_1,\dots,r_N\}$ is at least two and
equals the rank of $G$ as a free abelian group.
In this case,
the lengths are maximally dispersed in $\Real^*_+$.
Typically,
their multiplicity still grows exponentially,
but at a slower rate than in the nongeneric nonlattice case or the lattice case.

The geometry of a nonlattice self-similar string is much more elusive.
Indeed,
the only way to obtain detailed information about the volume of the tubular neighborhoods might be via the geometric zeta function,
the complex dimensions,
and the explicit formula~(\ref{E: Ve exp}).
Much like a real number can be approximated by a suitable sequence of rational numbers,
obtained by means of its continued fraction expansion (see~\cite{cxdss} and~\cite{HW,Schmidt}),
the set of lattice strings is dense in the space of all self-similar strings.
By an explicit Diophantine approximation procedure 
(see~[\citen{book},\,Thm.\,2.26,\,\citen{cxdss}]),
 we can approximate the nonlattice string $\L$ by a sequence
of lattice strings $\{\L^{(n)}\}_{n=1}^\infty$ with oscillatory
 period~$\oscp_n=2\pi q_n/\log r^{-1}$,
for larger and larger integers $q_n$.
Here,
$r=r_1$ is the first (i.e., largest) self-similarity ratio of $\L$.
Hence we can approximate,
in the sense of~[\citen{book},\,Def.\,2.5],
the set of complex dimensions $\D_\L$ of $\L$ by the corresponding sets of complex
 dimensions~$\D_n=\D_{\L^{(n)}}$.
In particular,
the complex dimensions of a nonlattice string exhibit a {\em quasiperiodic pattern,\/}
as explained in the following counterpart of Theorem~\ref{T: l}
(see~[\citen{book},\,Thms.\,2.13,\,6.20],
 along with~\cite{cxdss}):

\begin{theorem}[Nonlattice Case: Quasiperiodic Patterns]\label{T: nl}
Let\/ $\L$ be a nonlattice self-similar string.
Then $\L$ has infinitely many complex dimensions (with positive real part).
The density of the complex dimensions (counted with multiplicity) is
\begin{gather}\label{E: dens}
\#\{\omega\in\D_\L\colon0\leq\Im\omega\leq T\}\leq\frac{\log r_N^{-1}}{2\pi}T+O(1),\quad\text{ as }T\to\infty,
\end{gather}
with equality if\/ $\L$ has only one gap.
Further,
$D$ is the only complex dimension with real part equal to $D$,
and the other complex dimensions lie in a horizontally bounded strip~$\sigma_l\leq\Re\omega<D,$
where
$\sigma_l=\inf\{\Re\omega\colon\sum_{j=1}^Nr_j^\omega=1\}$.\footnote
{For a generic nonlattice string it can be shown that $\sigma=\sigma_l$ is the unique real
solution to the equation $1+r_1^\sigma+\dots+r_{N-m-1}^\sigma=mr_N^\sigma$,
where $m=\#\{j\colon r_j=r_N\}$ is the multiplicity of~$r_N$.
Moreover,
in most cases,
including the generic nonlattice case,
the infimum is not attained.
In fact,
it is possible that the infimum is only attained in the lattice case.\label{F: sigma_l}}

Each complex dimension $\omega_u$ of\/ $\L$ gives rise to a sequence of complex dimensions~$\omega$ of\/ $\L$,
close to the sequence~(\ref{E: line}).
Let\/ $r_1$ denote the largest scaling ratio of\/ $\L$.
Then the period of this sequence is $\oscp=2\pi/\log r_1^{-1}$.
Further,
 for certain multiples of this oscillatory period,
$\oscp_n=q_n\oscp=2\pi q_n/\log r_1^{-1}$,
where the sequence of integers~$q_n$ depends on the arithmetic nature of the scaling ratios
and is such that\/ $q_n\to\infty$ as~$n\to\infty$,
one obtains subsequences that lie arbitrarily close to the sequence~(\ref{E: line}),
with an oscillatory period of\/ $q_n\oscp$ instead of\/ $\oscp$.
\end{theorem}

We see that nonlattice complex dimensions `lie on infinitely many vertical lines' with larger and larger oscillatory period.
Thus we have an infinite spectrum of self-similarity in the nonlattice case.

The following theorem refines the information about the complex dimensions close to the sequence
$\{D+ki\oscp\colon k\in\Int\}$ when~$N=2$.
There are two scaling ratios~$r_1$ and $r_2=r_1^\alpha$,
for some irrational number $\alpha>1$,
and we write
\begin{gather}\label{E: f}
f(s)=1-r_1^s-r_2^s=1-r_1^s-r_1^{\alpha s},
\end{gather}
so that $f'(D)=(r_1^D+\alpha r_2^D)\log r_1^{-1}$.
See~[\citen{cxdss},\,Thm.\,4.3] for details and a proof.

\begin{theorem}[Dimensions above $D$]\label{T: Dline}
For $k\in\Int,$
let\/ $x=2\pi i(k\alpha -l),$
where $l$ is the integer nearest to $k\alpha$.
Then the complex dimension of\/ $\L$ close to $D+ik\oscp,$
 $\oscp=2\pi/\log r_1^{-1},$
is approximated by
\begin{gather}\label{E: omega}
\omega =D+ki\oscp-\frac{r_2^{D}}{f'(D)}\, x+\frac{(\log r_1)^2r_1^{D} r_2^{D}}{2{f'(D)}^3} \, x^2
+O(x^3).
\end{gather}
\end{theorem}

Note that the first correction term (of order $x$) is purely imaginary,
and the second correction is negative of order $x^2$ (since $x^2<0$).
Hence we see again that $\Re\omega<D$.
We refer to~[\citen{cxdss},\,Eq.\,(4.16)] for the 
general case when $N>2$.

This theorem allows us to obtain a dimension-free region,
recalled in Equation~(\ref{E: dimfree}) below.
Using the explicit formula~(\ref{E: Ve exp}) and this dimension free region,
we obtain,
by techniques explained in~\cite{pot},
that for every (small) $\delta >0$,
\begin{gather}\label{E: V leading order}
V(\varepsilon )
=\M\varepsilon^{1-D}\left( 1+O\left( \frac1{|\log\varepsilon |^{(N-1)/2-\delta}}\right)\right),\quad
\text{as }\varepsilon\to 0^+.
\end{gather}
The error estimate in~(\ref{E: V leading order}) is best possible when the system
 $\log r_j/\log r_1$ ($j=2,\dots,N$) is badly approximable in the sense of~\S\ref{S: dfree}.

It follows from~(\ref{E: V leading order}) and a computation of the residue in~(\ref{E: MDL})
that a nonlattice self-similar string is Minkowski measurable,
with a Minkowski content
\begin{gather}\label{E: M nl}
\M=\M(D;\L)=\frac{2^{1-D}L^D\left(g_1^D+\dots+g_K^D\right)}{D(1-D)\left(r_1^D\log r_1^{-1}+\dots+r_N^D\log r_N^{-1}\right)}.
\end{gather}

\begin{remark}
(a)
The statement of Theorem~\ref{T: nl} regarding the Minkowski measurability of nonlattice
 self-similar strings (i.e.,
of nonlattice self-similar sets in $\Real$) was first obtained independently in~\cite{mL93} and by Falconer in~\cite{kF95},
by means of the Renewal Theorem~[\citen{wF66},\,Ch.\,XI].
It has recently been extended to any (suitable) nonlattice self-similar set in $\Real^d$ ($d>1$)
by Gatzouras in~\cite{Gat},
thereby proving the nonlattice case of the geometric part of~[\citen{mL93},\,Conj.\,3,\,p.\,163] concerning self-similar drums.
(See also~[\citen{book},\,Conj.\,10.13,\,p.\,209] for an extension of part of that
 conjecture to the more general `tube formulas' obtained in the context of the theory of complex dimensions.)

(b)
Part of the results of~\cite{book} concerning the complex dimensions of self-similar strings was extended to random
self-similar strings (and other random fractal strings,
such as the zero set of Brownian motion) by Ben Hambly and the first author in~\cite{bHmL99}.
A rich variety of behaviors of the complex dimensions is encountered in this context.
\end{remark}

\subsection{Nonlattice Strings: Examples and Open Problems}
\label{S: nl ex conj}

In order to illustrate Theorem~\ref{T: nl} and to formulate some problems and a conjecture
 regarding the nonlattice case (generic and nongeneric),
we discuss three examples of nonlattice strings.
Our first example is a nongeneric nonlattice string.

\begin{example}[A Nongeneric Nonlattice String]\label{E: nong ex}
Consider the self-similar string with scaling ratios
$r_1=1/2$,
$r_2=1/4$,
$r_3=2^{-1-\sqrt{2}}$ and one gap $g=1/4-r_3$.
Figure~\ref{F: stages} gives six approximations to its complex dimensions
(along with an associated graph of the density of their real parts in the lower diagrams,
to be explained below)
corresponding to the successive approximations to $\sqrt2$ (obtained by the continued fraction expansion~\cite{HW}):
$\frac{p_n}{q_n}=\frac32, \frac75, \frac{17}{12}, \frac{41}{29}, \frac{99}{70}, \frac{239}{169},\dots\to\sqrt{2}$.
For example,
the first approximation gives the lattice string with scaling ratios $\tilde r_1=2^{-1}$,
$\tilde r_2=2^{-2}$,
$\tilde r_3=2^{-5/2}$,
the complex dimensions of which are the solutions to the equation~$z^2+z^4+z^5=1$,
$2^{-\omega/2}=z$.
The oscillatory period of this lattice approximation is $\oscp=4\pi/\log 2$.
One sees the development of a quasiperiodic pattern:
the complex dimensions of the nonlattice string are well approximated by those of a lattice string for a certain {\em finite\/}
number of periods of the lattice approximation.
Then that periodic pattern gradually disappears,
and a new periodic pattern,
approximated by the next lattice approximation,
emerges.
\end{example}

\begin{figure}[tbh]
\begin{picture}(360,335)(0,0)
\put(0,0){\put(66.8,180){$\frac\oscp2$}
	\epsfxsize 4.2cm\epsfbox{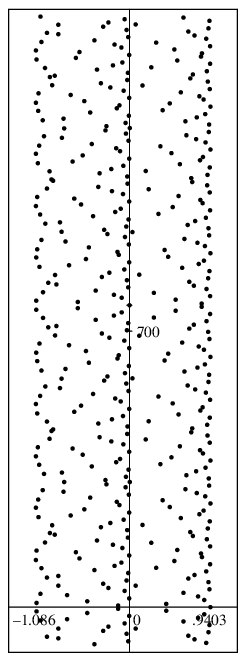}}
\put(120,0){\put(66.8,180){$\frac\oscp2$}
	\epsfxsize 4.2cm\epsfbox{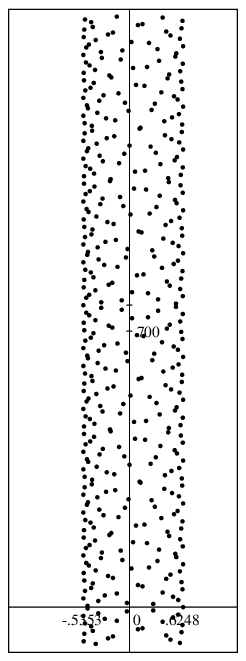}}
\put(240,0){\put(66.8,181){$\scriptstyle\oscp$}
	\epsfxsize 4.2cm\epsfbox{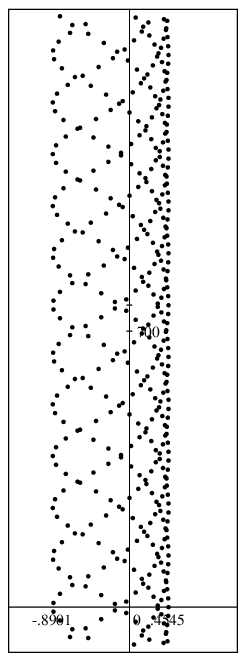}}
\end{picture}
\caption{An impression of the quasiperiodic pattern of complex dimensions on a larger scale,
for the self-similar string of Example~\ref{E: nong ex} (left)
and those of Example~\ref{E: gen ex} (center and right).}
\label{F: large scale}\end{figure}

It takes fairly large approximations to see this quasiperiodic pattern.
For example,
the complex dimensions of the three nonlattice strings of Examples~\ref{E: nong ex}
and~\ref{E: gen ex} are respectively approximated to within a distance $\frac1{10}$ by the points
 in the three diagrams of Figure~\ref{F: large scale} for only about two periods $\oscp$ of the
respective lattice strings (four periods for the last one,
since $\oscp$ is half the size,
as explained in Example~\ref{E: gen ex}).
This is found by adapting the proof of~[\citen{cxdss},\,Thm.\,3.6] to these strings,
since a direct application of this theorem would give only half a period of good approximation.
The number of periods (of the lattice string) for which the approximation is good does grow
 linearly in the denominator of the approximation.
Note that the period itself also grows like this denominator.

\begin{figure}[tbh]
\begin{picture}(360,100)(0,-20)
\put(0,0){\put(60,-12){\makebox(0,0){$r_1=\frac12$, $r_2=\frac14$, $r_3=2^{-1-\sqrt{2}}$}}
	\epsfxsize 4.1cm\epsfbox{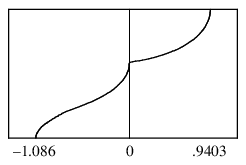}}
\put(120,0){\put(60,-12){\makebox(0,0){$r_1=\frac12$, $r_2=2^{-1-\sqrt{2}}$}}
	\epsfxsize 4.1cm\epsfbox{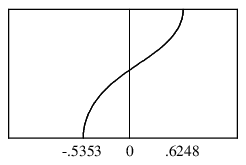}}
\put(240,0){\put(60,-12){\makebox(0,0){$r_1=\frac14$, $r_2=2^{-1-\sqrt{2}}$}}
	\epsfxsize 4.1cm\epsfbox{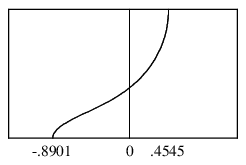}}
\end{picture}
\caption{Comparison of the densities of the real parts for the nongeneric nonlattice string
of Example~\ref{E: nong ex} (left) and the two generic nonlattice strings of
 Example~\ref{E: gen ex} (center and right).}
\label{F: density}\end{figure}

The left diagram in Figure~\ref{F: large scale} gives an impression on a larger scale of the
 quasiperiodic behavior of the complex dimensions of the nongeneric nonlattice string of
Example~\ref{E: nong ex}.
The other two diagrams in this figure give the complex dimensions of the two {\em generic\/}
 nonlattice strings of Example~\ref{E: gen ex} below.
One sees that the complex dimensions in the left diagram are much denser to the left of vanishing real part.
One sees this even more clearly in Figure~\ref{F: density},
where the cumulative density of the real parts is graphed for these three self-similar strings.
We have no explanation for this apparent `phase transition'.
We formulate this question as a problem
(see Theorem~\ref{T: nl} and footnote~\ref{F: sigma_l} for the definition of $\sigma_l$):

\begin{problem}[Transition in the Nongeneric Nonlattice Case]
A nongeneric nonlattice string has a vertical line of transition inside the strip
 \mbox{$\sigma_l\leq\Re s\leq D$},
to the left of which the density of the real parts is infinitely higher than to the right.
Such a transition does not occur for generic nonlattice strings.
\end{problem}

Thus,
for the nongeneric nonlattice string of Example~\ref{E: nong ex},
this transition occurs at \mbox{$\Re s=0$},
as indicated by the corner of the density graph at this point,
and the vertical part of the graph to the left of~$\Re s=0$.
Moreover,
from other numerical evidence,
it seems that this line of transition often occurs at \mbox{$\Re s=0$}.

\begin{example}[Two Generic Nonlattice Strings]\label{E: gen ex}
We also include
in Figures~\ref{F: large scale} and~\ref{F: density}
 the analogous diagrams for two generic nonlattice strings.
These are respectively the nonlattice strings with two scaling ratios $r_1=1/2$ and
 $r_2=2^{-1-\sqrt{2}}$ (the middle diagrams in both figures) and $r_1=1/4$ and
 $r_2=2^{-1-\sqrt{2}}$ 
(the right diagrams in both figures).
Each of these strings has a single gap of length respectively $g_1=1/2-2^{-1-\sqrt{2}}$ and $g_1=3/4-2^{-1-\sqrt{2}}$.
Note that for the second string,
with $r_1=1/4$ and $r_2=2^{-1-\sqrt{2}}$,
the approximation $1+\sqrt{2}\approx\frac{408}{169}$ leads to the equation
\begin{gather*}
z^{2\cdot169}+z^{408}=1,\qquad 2^{-\omega/169}=z.
\end{gather*}
Since this is an equation in $z^2$,
the oscillatory period of the lattice equation,
namely,
\mbox{$169\cdot\pi/\log 2$},
 is only slightly larger than that of the lattice string corresponding to the previous approximation,
$1+\sqrt{2}\approx\frac{169}{70}$,
which is $70\cdot2\pi/\log 2$.
\end{example}

We see in the Figures~\ref{F: large scale} and~\ref{F: density} that the complex dimensions of the
 nongeneric and generic nonlattice strings of 
Examples~\ref{E: nong ex} and~\ref{E: gen ex} are much denser near the boundaries $\Re s=\sigma_l$
and $\Re s=D$ of the critical strip.
Indeed,
in~\cite{cxdss},
we have shown that given any nonlattice string with two scaling ratios,
there exists an explicit positive constant
 $\refstepcounter{constant}\label{c: dens}C_{\ref{c: dens}}$
(expressed in terms of $r_1$ and $r_2$),
such that the density graph of the real parts is approximated by the graph
 $y=1-C_{\ref{c: dens}}\sqrt{D-x}$,
 for $x\leq D$ in a small neighborhood to the left of $D$.
More generally,
if the rank of $G$ is~$v$ (so $v\leq M$,
and $v=M$ in the generic nonlattice case,
see \S\ref{S: l}),
then the density graph is approximated by the graph 
$\refstepcounter{constant}\label{c: densM}y=1-C_{\ref{c: densM}}(D-x)^{(v-1)/2}$.
Hence,
if there are three multiplicatively independent scaling ratios ($v=M=3$),
then the density at the boundary is comparable to the average density,
whereas for $v\geq4$,
the complex dimensions thin out near the boundary of the `critical strip'
 $\sigma_l\leq\Re\omega\leq D$.
Indeed,
numerical evidence suggests the following problem,
further evidence for which will be provided in~\S\ref{S: dfree},
in connection with the properties of Diophantine approximation of the ratios $\log r_j/\log r_1$.
We refer to~[\citen{cxdss},\,Rem.\,4.14] for more details.
Also see~\S\ref{S: dfree} for some information about bad approximability.

\begin{problem}[Complex Dimensions Close to One Line]\label{P: single line}
The complex dimensions of a nonlattice string with a large number of different scaling ratios
 are concentrated around a single line within its critical strip.
That is,
the graph of the density of the real parts of the complex dimensions has a narrow steep part,
and two flat parts to the left and the right.
Moreover,
if this phenomenon occurs,
it may be most apparent when the scaling factors~$r_j$ are such that the ratios
 $\log r_j/\log r_1$ ($j=2,\dots,N$) are badly approximable.
\end{problem}

We have shown in~[\citen{cxdss},\,Thm.\,8.1] that for any nonlattice string,
 the real parts of the complex dimensions form a set without isolated points.
In~[\citen{book},\,Ch.\,10],
we proposed as a definition of fractality
the existence of at least one nonreal complex dimension with positive real part.
Since nonreal complex dimensions give oscillatory terms in the explicit formula for the tubular neighborhoods,
that is,
terms that are `self-similar' under a (multiplicative) scaling,
this theorem can be phrased as saying that
a nonlattice string is fractal in a perfect set of dimensions.
The following conjecture even says that a nonlattice string is fractal in every dimension
 in~$[\sigma_l,D]$
(see Theorem~\ref{T: nl} for the definition of $\sigma_l$).

\begin{conjecture}[Density of the Real Parts]\label{C: dens}
Let $\L$ be a nonlattice string.
Then the real parts of its complex dimensions form a set that is dense in the connected interval $[\sigma_l,D]$.
\end{conjecture}

We sketch the idea for a possible proof of this conjecture when $N=2$.
By Equation~(\ref{E: omega}),
one can find a complex dimension with real part arbitrarily close to $D$.
Moreover,
for every complex dimension $\omega$ there exists an analogue of Equation~(\ref{E: omega}),
yielding a sequence of complex dimensions close to the line $\omega+ki\oscp$ ($k\in\Int$).
The main correction term,
of order $x$,
will in general not be purely imaginary,
as it is for $\omega=D$,
and for suitable values of $k$,
the real part of $x$ will be negative.
Continuing in this manner,
given any $\varepsilon>0$,
one can find a sequence of complex dimensions,
starting at $D$,
the real parts of which form a decreasing sequence of numbers that differ by less
 than~$\varepsilon$.
This would prove the density of real parts in an interval.
It remains to establish a suitable version of Equation~(\ref{E: omega})
 (and its analogues for~$N>2$),
for every complex dimension,
and to analyze it.
Of course,
the argument that the real part of $x$ is negative for suitable values of $k\in\Int$ breaks down when one reaches the left boundary
of the `critical strip'.
There,
$x$ will become purely imaginary (like at $D$),
and the second correction term will become positive.
Note,
 however,
that there is no complex dimension with real part equal to $\sigma_l$.

\section{Dimension-Free Regions}
\label{S: dfree}

\begin{definition}
An open domain in the complex plane is a {\em dimension-free region\/} for the string
$\L$ if it contains the line $\Re s=D$ and the only pole of $\zeta_\L$ in
that domain is $s=D$.
\end{definition}

In the situation of Theorem~\ref{T: Dline},
when $N=2$,
the dimension-free region depends on the approximability by rationals of $\alpha$,
where $\alpha=\log r_2/\log r_1$ as in Equation~(\ref{E: f}).
Assume that $\alpha$ is {\em badly approximable\/}:
$\refstepcounter{constant}\label{C: cf}|q\alpha-p|\geq C_{\ref{C: cf}}/q$ for all integers $p$
and~$q$ (see,
e.g.,~\cite{Dodson}).\footnote
{Equivalently,
$\alpha$ has bounded partial quotients in its continued fraction expansion.}
Put $\refstepcounter{constant}\label{C: b}C_{\ref{C: b}}=\pi^4(r_1r_2)^D/(2{f'(D)}^3)$,
where~$f$ is the function defined in~(\ref{E: f}).
Then~\mbox{$
\left\{\sigma +it\in\Com\colon\sigma >D-C_{\ref{C: b}}C_{\ref{C: cf}}^{2}t^{-2}\right\}
$} is a dimension-free region for~$\L$.
More generally,
let\/ $g\colon\Real_+\rightarrow (0,1]$ be a nonincreasing function that measures the
 approximability of $\alpha$ by rationals:
\mbox{$|q\alpha-p|\geq g(q)/q$} for all integers $p$ and~$q$.
Write $w_1=\log r_1^{-1}$.
Then $\L$ has a dimension-free region of the form
\begin{gather}
\left\{\sigma +it\in\Com\colon\sigma >D-C_{\ref{C: b}}g^{2}(w_1t/2\pi)t^{-2}\right\} .
\end{gather}

In the general case,
when $N\geq2$,
the dimension-free region is best (widest) if the ratios $\log r_j/\log r_1$ ($j=2,\dots,N$) are badly approximable.
In that case,
$\L$ has a dimension-free region of the form
 (see~[\citen{pot}--\ref{cxdssc}])\refstepcounter{constant}\label{C: bN}
\begin{gather}\label{E: dimfree}
\bigl\{ \sigma +it\in\Com\colon\sigma> D-C_{\ref{C: bN}}t^{-2/(N-1)}\bigr\} ,
\end{gather}
where the constant $C_{\ref{C: bN}}$ depends only on $r_1,\dots ,r_N.$
On the other hand,
if those ratios are better approximable,
then the dimension-free region is worse,
in the sense that there are complex dimensions that lie very close to the vertical line $\Re s=D$.

From these results it is also clear that in general,
depending on the properties of simultaneous Diophantine approximation of the scaling ratios,
nonlattice strings with a larger number of scaling ratios have a wider (and hence better) dimension-free region.
This reinforces the belief in Problem~\ref{P: single line}.
For more information regarding the results of \S\ref{S: dfree} and their proofs,
as well as for related results also based in part on Diophantine approximation,
we refer to~[\citen{pot}\,\&\,\ref{cxdssc},\,\S7].

\section{Self-Similar Dynamical Systems}
\label{S: dyn}

We refer the interested reader
 to~[\citen{sL89},\,\citen{pot}--\ref{cxdssc},\,\citen{wPmP83}--\ref{wPmP90c}],
and the relevant references therein,
for more information on the dynamical systems considered here.

Regarding earlier results on the Prime Orbit Theorem in a related context,
we mention,
in particular,
the early works of Huber~\cite{Hu},
Sinai~\cite{Si},
and Margulis~\cite{Mar}
(see the historical note in~[\citen{BedKS},\,p.\,154]) and the more recent works
by Parry and Pollicott~[\citen{wPmP83}\,\&\,\ref{wPmP90c},\,Ch.\,6] as well as by
 Lalley~\cite{sL89}.

For simplicity and because this is the situation that is best understood so far,
we consider here the dynamical counterpart of self-similar strings with a
 single gap~(
$K=1$).
It would be interesting to extend our formalism to include the multiple gap case.
This can possibly be done within the setting of~\cite{bHmL99} where `Galton-Watson'
random trees are used to study the zeta functions of random strings.
\medskip

Let $N\geq 0$ be an integer and let $\Omega =\{ 1,\dots ,N\}^\Nat$ be the space of sequences over the alphabet
$\{ 1,\dots ,N\}$.
Let $\weight\colon\Omega\rightarrow (0,\infty ]$ be a function,
called the {\em weight.\/}
On~$\Omega$,
we have the left shift $\tau$,
given on a sequence $z=(z_n)$ by $(\tau z)_n=z_{n+1}$.
We define the {\em suspended flow\/} $\Flow_\weight$ on the space $[0,\infty )\times\Omega$ as the following
continuous-time dynamical system
(see~[\citen{pot},\,Fig.\,1,\,p.\,117]):
\begin{gather}
\Flow_\weight (t,z)=\begin{cases}
(t,z)			&\mbox{if }0\leq t<\weight(z),\\
\Flow_\weight (t-\weight(z),\tau z)	&\mbox{if }t\geq \weight(z).
\end{cases}
\end{gather}

Given a finite sequence $z_1,z_2,\dots ,z_l$,
we let
$z=z_1,z_2,\dots ,z_l,z_1,z_2,\dots ,z_l,\dots$
be the corresponding periodic sequence,
and
$\orb=\{z,\tau z,\tau^2z,\dots\}$
the associated finite orbit of $\tau$,
of length $\#\orb$.
Thus the length of $\orb$ is a divisor of $l$.
The {\em total weight\/} of $\orb$ is
$\tweight (\orb)=\sum_{z\in\orb}\weight(z)$.
We define the {\em dynamical zeta function\/}
 (or Bowen--Ruelle zeta function~\cite{rB73,dR76,wPmP90})
 of $\Flow_\weight$ by the Euler-type product
\begin{gather}
\zeta_\weight(s)=\prod_\orb\frac1{1-e^{-\tweight(\orb)s}},
\end{gather}
where $\orb$ runs over all finite orbits of $\tau$.

The function
$\psi_\weight (x)=\sum_{k\tweight (\orb )\leq\log x}\tweight (\orb )$
 is the analogue of the function $\psi(x)=\sum_{p^k\leq x}\log p$ in number theory,
which counts the prime powers with a weight $\log p$.\footnote
{In the literature,
$\psi$ is traditionally referred to as the `von Mangoldt function' and is often used,
in particular,
to prove the classical Prime Number Theorem;
see,
e.g.,
[\citen{In,Pat,Ti}].}
It counts the periodic orbits of $\Flow_\weight$ and their multiples by their total weight.
It may be referred to as the {\em weighted prime orbit counting function\/} and is related to 
the dynamical zeta function by
$-\zeta_\weight'(s)/\zeta_\weight(s)=\int_0^\infty x^{-s}d\psi_\weight(x)$,
for $\Re s>D$.
The poles of the logarithmic derivative $\zeta_\weight'/\zeta_\weight$ are called the 
{\em dynamical complex dimensions\/} of $\Flow_\weight$.
Thus the dynamical complex dimensions are the poles {\em and the zeros\/} of $\zeta_\weight$,
counted {\em without\/} multiplicity.
(See~\cite{pot} for more information.)
\medskip

The flow $\Flow_\weight$ is said to be {\em self-similar\/} if\/ $N\geq 2$ and the weight function\/ $\weight$ depends only
on the first letter of the sequence at which it is evaluated.
We then define the scaling ratios of $\Flow_\weight$ by \mbox{$r_j=\exp(-\weight(j,j,\dots ))$}
for $j=1,\dots ,N$.
Much as in~\S\ref{S: l nl},
the flow is said to be {\em lattice} if the group $G$ defined in~\S\ref{S: l} has rank~$1$,
and {\em nonlattice} otherwise.
We associate with $\Flow_\weight$ a self-similar fractal string with these scaling ratios
and one gap $g=1-r_1-\dots-r_N$ and total length such that~\mbox{$gL=1$}
(we need to require that $r_1+\dots+r_N<1$).
Then we have~\mbox{$\zeta_\weight(s)=\zeta_\L(s)$}
(see~\S\ref{S: ss} and~[\citen{pot},\,Thm.\,2.10]).
Since $\zeta_\weight$ has no zeros in this case,
 the dynamical complex dimensions of a self-similar flow coincide with the geometric complex
 dimensions of the corresponding self-similar string,
but they are counted without multiplicity.

If\/ $D$ is the only complex dimension on the line $\Re s=D$,
then we deduce (from the explicit formula for the prime orbit counting function)
 a {\em Prime Orbit Theorem\/} for~$\Flow_\weight$:
\mbox{$\psi_\weight (x)=x^D/D+o\left( x^D\right)$},
as $x\to\infty$.
This is the analogue of the famous Prime Number Theorem,
which can be stated as $\psi(x)=x+o(x)$.
In the nonlattice case,
using the dimension-free region~(\ref{E: dimfree}) and the techniques explained in~\cite{pot},
we even obtain for every (small)~$\delta>0$ the following
 {\em Prime Orbit Theorem with Error Term\/}:
\begin{gather}\label{E: psi leading order}
\psi_\weight (x)=\frac{x^D}{D}\left(1+O\left(\frac1{(\log x)^{(N-1)/2-\delta}}\right)\right) ,\quad
\text{as }x\to\infty.
\end{gather}
This estimate is best possible if the system
$\{\log r_j/\log r_1\colon j=2,\dots,N\}$ is badly approximable 
(as explained for $N=2$ in~\S\ref{S: dfree}).
We note that in~\cite{pot},
we also obtain more precise explicit formulas for the counting function $\psi_\weight$.

\section{The Complex Dimensions as the Spectrum of Shifts}
\label{S: shifts}

In the simplest case of the Cantor string,
discussed in the introduction and in Example~\ref{E: mCS},
it is easy to see that a shifted copy,
 $\CS+x$,
 overlaps the Cantor string only for shifts over~$2\cdot3^{-n}$,
$n=0,1,\dots$,
or suitable combinations of such shifts:
$x=\sum_{n=0}^\infty a_n3^{-n}$,
 $a_n=\pm2$.
For all other values of the real number $x$,
there is no overlap:
$\CS\cap(\CS+x)=\emptyset$.
Moreover,
the Minkowski dimension (or Hausdorff dimension) of the intersection is always $D=\log_32$ and its Minkowski content (or Hausdorff measure)
 is $2^k$ times smaller than that of the Cantor set itself,
where $k$ is the number of nonzero digits of $x$ in the above representation.
This information suffices to reconstruct the Cantor set and its complex dimensions.

In general,
let $F$ be a self-similar fractal subset of the real line (see Remark~\ref{R: ss boundary}).
Let $I_\varepsilon$ denote the interval $(-\varepsilon,\varepsilon)$.
To measure the overlap of $F$ with a shifted copy of it,
we consider the function\footnote
{See footnote~\ref{F: set sum} for the notation $F+I_\varepsilon$.}
$
f(\varepsilon,x)=\left|(F+I_\varepsilon)\cap(F+I_\varepsilon+x)\right|$,
 for $x\in\Real$.
We then construct the `dimension-like function'
$$
D(x)=\inf\{d\geq0\colon f(\varepsilon,x)=O\left(\varepsilon^{1-d}\right)\text{ as }\varepsilon\to0^+\}
$$
and the `upper Minkowski content-like function'
$$
\M^*(x)=\limsup_{\varepsilon\to0^+}f(\varepsilon,x)\varepsilon^{D(x)-1}.
$$
We expect that these functions exhibit the following behavior.
In the lattice case,~$D(x)>0$ for a bounded discrete set of values.
The function $\M^*(x)$ is discontinuous at each of these values,
and continuous and vanishing on the complement.
In the nonlattice case,
$D(x)>0$ for every $x$ in a countable dense subset of a compact connected interval and $\M^*(x)$ is continuous.
The complex dimensions of $F$ can be recovered from the functions $D(x)$ and $\M^*(x)$.
We hope to develop these ideas in subsequent work.

\section{A Cohomological Interpretation of the Complex Dimensions}
\label{S: coho}

We have suggested in~[\citen{book},\,\S10.5] that there should exist a notion of
 `complex cohomology' that could be associated to both 
arithmetic and self-similar geometries and have proposed a possible dictionary towards it,
building in part upon the analogy between lattice strings (or more generally,
lattice self-similar geometries) and curves (or varieties) over finite
 fields~[\citen{aW41}--\ref{aW52c},\,\citen{aPiS95},\,I].
We have proposed,
in particular,
to construct a `cohomology theory',
which to each (dynamical) complex dimension $\omega$ of such a geometry,
would associate a nontrivial cohomology group~$H^\omega$ (with coefficients in the field of complex numbers).
In general,
$H^\omega$ would be expected to be an infinite dimensional Hilbert space.
This should be the case,
for example,
for the cohomology spaces associated with a nonlattice string (or more generally,
with a nonlattice self-similar geometry).
A possible way of developing such a theory was proposed at the very end of~[\citen{book},\,\S10.5].
It relies on a suitable equivalence relation between the Dirichlet series (or `zeta functions') attached to these geometries;
see~[\citen{book},\,p.\,220].

\subsection{Fractal Geometries and Finite Geometries}

The geometric zeta function of a fractal string can have complex dimensions of higher multiplicity,
like in Example~\ref{E: not simple} above.
Hence $\zeta_\L(s)$ should be compared to the zeta function of a variety over a finite field,
and not for example to its logarithmic derivative.
This is confirmed by the fact that the residue at a complex dimension is quite arbitrary,
and not an integer in general.
Another confirmation is that the logarithmic derivative of the geometric zeta function is the
logarithmic derivative of the dynamical zeta function
(as shown in~\cite{pot} and recalled in~\S\ref{S: dyn} above),
which is the generating function of the counting function of the periodic orbits of a
 dynamical system.
This corresponds to the logarithmic derivative of the zeta function of a variety,
which is the generating function of the counting function of the periodic orbits of the Frobenius automorphism,
as explained in \S\ref{S: ss finite}.
On the other hand,
for the simplest self-similar strings (i.e.,
those without gaps~\cite{book}),
the geometric zeta function has only poles and no zeros.
This would mean that there is only `even dimensional cohomology' in that case,
but the correct interpretation is as yet unclear.
\medskip

We complete the dictionary proposed in~[\citen{book},\,\S10.5] with the following three tables.
The first table summarizes the analogies between the corresponding zeta functions.
\medskip

\newcounter{tbl}\setcounter{tbl}{0}
\begin{center}
\begin{tabular}{ll}
\em self-similar geometries	&\em `finite' geometries\\
\hline
lattice string			&variety $V$ over the field $\Fin_q$\\
nonlattice string $\L$		&infinite dimensional variety\\
Cantor string			&affine rational variety\\
geometric zeta function $\zeta_\L$	&zeta function $\zeta_V$\\
\ind counts lengths		&\ind counts divisors\\
lattice case is periodic	& the zeta function is periodic\\
period $\oscp$ is typically large	&period is typically small\\
nonlattice: $\oscp\to\infty$	&period $\to0$ if $\#\Fin_q\to\infty$\\
\ind (i.e., `characteristic' $\to1$)\\
residue at $D$ (Minkowski	&residue at $n$ (dimension)\\
\ind dimension) gives $\M(D;\L)$	&\ind gives the class number\\
\end{tabular}\\[4mm]
\refstepcounter{tbl}\label{T: zeta}
\begin{minipage}{105mm}
{\sc Table} \ref{T: zeta}. Self-similar fractal geometries vs.\ varieties over finite fields:
Analogies between the zeta functions.
\end{minipage}
\end{center}
\medskip

We explain the sense in which $\oscp\to\infty$ for nonlattice strings,
and the connection with the characteristic.
Recall that the oscillatory period of a lattice string with scaling ratios
$r_j=r^{k_j}$,
$j=1,\dots,N$ (for positive integers $k_1,\dots,k_N$ without common factor and $r\in(0,1)$) is 
$\oscp=\frac{2\pi}{\log r^{-1}}$.
It is the period of its geometric zeta function:
$\zeta_\L(s)=\zeta_\L(s+i\oscp)$,
for $s\in\Com$.
On the other hand,
if $\L$ is a nonlattice string,
and $r_j\approx r^{k_j}$ gives a lattice approximation,
then $\zeta_\L(s)$ and $\zeta_\L(s+in\oscp)$ are close,
for~$n$ not too large.
The next (better) lattice approximation to $\L$ has larger values for the integers $k_j$ ($j=1,\dots,N$) and a value of $r$ that is closer to $1$.
Hence,
 its oscillatory period is larger (see \S\ref{S: l nl} for more details).
In that sense,
$\oscp\to\infty$ for nonlattice strings,
as $\oscp$ runs through the values of the oscillatory periods of lattice approximations.
For a variety (of dimension $n$) over the finite field $\Fin_q$,
 on the other hand,
the zeta function is periodic with period
$\oscp=\frac{2\pi}{\log q}$.
Again,
the zeta function of $V$ is periodic:
$\zeta_V(s)=\zeta_V(s+i\oscp)$,
for $s\in\Com$.
Here,
$q$ is a power of a prime number,
the characteristic of the finite field.
Hence the period is small typically.
However,
the limit as~$q\downarrow1$ (i.e.,
if the characteristic of $\Fin_q$ tended to $1$,
if this were possible),
corresponds to the limit $r\uparrow1$ that we see for nonlattice strings.

The next table summarizes the expected properties of a cohomology theory for self-similar fractal strings (and sets).
\medskip

\begin{center}
\begin{tabular}{ll}
\em self-similar geometries	&\em `finite' geometries\\
\hline
poles from the scaling ratios	&poles from even cohomology\\
zeros from the gaps		&zeros from odd cohomology\\
lattice case is periodic	&zeta function is periodic\\
nonlattice cohomology is	& the cohomology collapses to\\
\ind infinite dimensional	&\ind a finite-dimensional one\\
\end{tabular}\nopagebreak\\[4mm]\nopagebreak
\refstepcounter{tbl}\label{T: coho}\nopagebreak
\begin{minipage}{105mm}
{\sc Table} \ref{T: coho}. Self-similar fractal geometries vs.\ varieties over finite fields:
Cohomological aspects.
\end{minipage}
\end{center}
\medskip

Finally,
Table~\ref{T: dyn} presents the dynamical analogies between self-similar flows and the Frobenius flow of a variety.
\medskip

\begin{center}
\begin{tabular}{ll}
\em self-similar geometries	&\em `finite' geometries\\
\hline
dynamical flow			&Frobenius flow\\
dynamical zeta function $-\zeta_\L'/\zeta_\L$	& $-\zeta_V'/\zeta_V$ counts the Frobenius\\
\ind counts the closed orbits	&\ind or Galois orbits of points\\
Euler-type product connects	&Euler product connects \\
\ind orbits with lengths	&\ind orbits with divisors\\
$r^\omega$ (poles) are solutions	&$q^\omega$ (zeros and poles) are\\
\ind to~(\ref{E: eq})		&\ind eigenvalues of Frobenius\\
number of lines	is $k_N=$ degree	&number of lines is $2n+1$,\\
\ind  of (\ref{E: eq})		&\ind $n=\dim V$\\
\end{tabular}\\[4mm]
\refstepcounter{tbl}\label{T: dyn}
\begin{minipage}{105mm}
{\sc Table} \ref{T: dyn}. Self-similar fractal geometries vs.\ varieties over finite fields:
Dynamical aspects.
\end{minipage}
\end{center}

\begin{remark}\label{R: nl gaps}
One aspect of Table~\ref{T: coho} that we have not developed in this paper is that the dichotomy
 lattice vs.\ nonlattice also exists for the sequence of gaps.
A string may be a lattice string both for the  scaling ratios and the gaps
(and hence,
both the zeros and poles of $\zeta_\L$ are found by solving a polynomial equation,
by Equation~(\ref{E: zeta ss})),
but be nonlattice when the scaling ratios and gaps are considered together;
that is,
the numerator and denominator of $\zeta_\L$ have uncommensurable oscillatory periods.
\end{remark}

\end{document}